\pdfoutput=1

\documentclass[12]{amsart}

\usepackage{amsmath}
\usepackage{amsfonts}
\usepackage{rotating}
\usepackage{bbm}
\usepackage[all,cmtip]{xy}
\usepackage{subcaption}
\usepackage{enumitem}
\usepackage{longtable}
\usepackage{url,booktabs}
\usepackage{array}
\usepackage{ragged2e}
\newcolumntype{P}[1]{>{\RaggedRight\hspace{0pt}}p{#1}}

\usepackage{amssymb}
\usepackage[latin1]{inputenc}
\usepackage{graphicx}
\usepackage{dsfont}
\usepackage{color}
\usepackage{hyperref}
\usepackage{makecell}
\renewcommand{\phi}{\varphi}

\newcommand{\R}{{\mathbb{R}}}

\allowdisplaybreaks

\advance\hoffset by -0.5 truecm
\usepackage{amssymb}
\newtheorem{Theorem}{Theorem}[section]

\newtheorem{Remark}[Theorem]{Remark}

\def\R{{\mathbb R}}

\def \0{\lambda_{0}}

\begin{document}

\title[The Euler problem of two fixed centers]{Homoclinic orbits in the Euler problem of two fixed centers}

\author{Seongchan Kim}
 \address{Universit\"at Augsburg, Universit\"atsstrasse 14, D-86159 Augsburg, Germany}
 \email {seongchan.kim@math.uni-augsburg.de}
\date{\today}



\begin{abstract}

We give a complete description of the shapes and the behavior of all homoclinic orbits in the Euler problem of two fixed centers.

\end{abstract}

\maketitle
\setcounter{tocdepth}{1}
\tableofcontents

\section{Introduction}

A \textit{homoclinic orbit}, whose concept was introduced by Poincar\'e, is an orbit which is asymptotic to an unstable periodic orbit in both forward and backward time. In view of the Poincar\'e section map, an unstable periodic orbit corresponds to a hyperbolic fixed point of this map and a homoclinic orbit is represented by a point whose orbit is doubly  asymptotic to the  hyperbolic fixed point.

In this paper, we study homoclinic orbits in the Euler problem of two fixed centers. This problem  describes the motion of a point-mass under the influence of a Newtonian potential with two fixed attracting points. We refer to the two attracting bodies  as  the \textit{Earth} and  \textit{Moon} and the point-mass as the \textit{satellite}. This problem  can be obtained from the planar circular restricted three-body problem by ignoring the centrifugal and Coriolis terms. That this system is  integrable  was discovered by Euler in  1760.

We denote by $\mu \in (0,1)$ the mass ratio of the two massive bodies and locate the Earth and  Moon at ${E}=(-{1}/{2},0)$ and ${M}=({1}/{2},0)$, respectively. The  describing Hamiltonian $H : T^*(  \R^2 \setminus \left \{ E,M \right\}  )\cong ( \R^2 \setminus \left \{ {E, M} \right \} ) \times \R^2 \rightarrow \R$ is given by 
$$
H(q,p) = \frac{1}{2}|p|^2 - \frac{ 1-\mu}{|q-{E}|} - \frac{\mu}{|q-{M}|}.
$$
Without loss of generality we may assume that $\mu \leq {1}/{2}$, i.e., the Earth is stronger.  The Hamiltonian has a unique critical point $L = (l,0,0,0)$, $l \in (-{1}/{2},{1}/{2})$, of Morse index  1. It corresponds to the saddle point of the potential.

Given an energy level $H=c$, we define the \textit{Hill's region} to be 
$$\mathcal{K}_c := \pi ( H^{-1}(c) ) \subset \R^2 \setminus \left \{ {E, M} \right \},$$
 where $ \pi : ( \R^2 \setminus \left \{ {E, M} \right \} ) \times \R^2 \rightarrow \R^2 \setminus \left \{ {E, M} \right \}$ is the footpoint projection.  In what follows, we only consider negative energies so that     Hill's regions are   bounded, namely the satellite is confined to a bounded region in  the configuration space $ \R^2 \setminus \left \{ {E, M} \right \}$.  We denote by $c_J $ the critical energy level $H(L)= -1-2\sqrt{ \mu(1-\mu)}$.  For $c < c_J$, the Hill's region consists of two connected components, where one of them is a neighborhood of the Earth and the other is a neighborhood of the Moon. We denote them by $\mathcal{K}_c^E$ and $\mathcal{K}_c^M$, respectively.  For $c > c_J$, these two components become connected. The satellite then can move from  the Earth to the  Moon and vice versa. There are  two further distinguished energy levels $c_e$ and $c_h$ at which the Liouville foliation changes, see Figure \ref{diagram}. Note that $c_J < c_e < c_h < 0$ if $\mu \neq 1/2$ and $c_J < c_e < c_h=0$ if $\mu=1/2$.  

  \begin{figure}[t]
  \centering
  \includegraphics[width=0.85\linewidth]{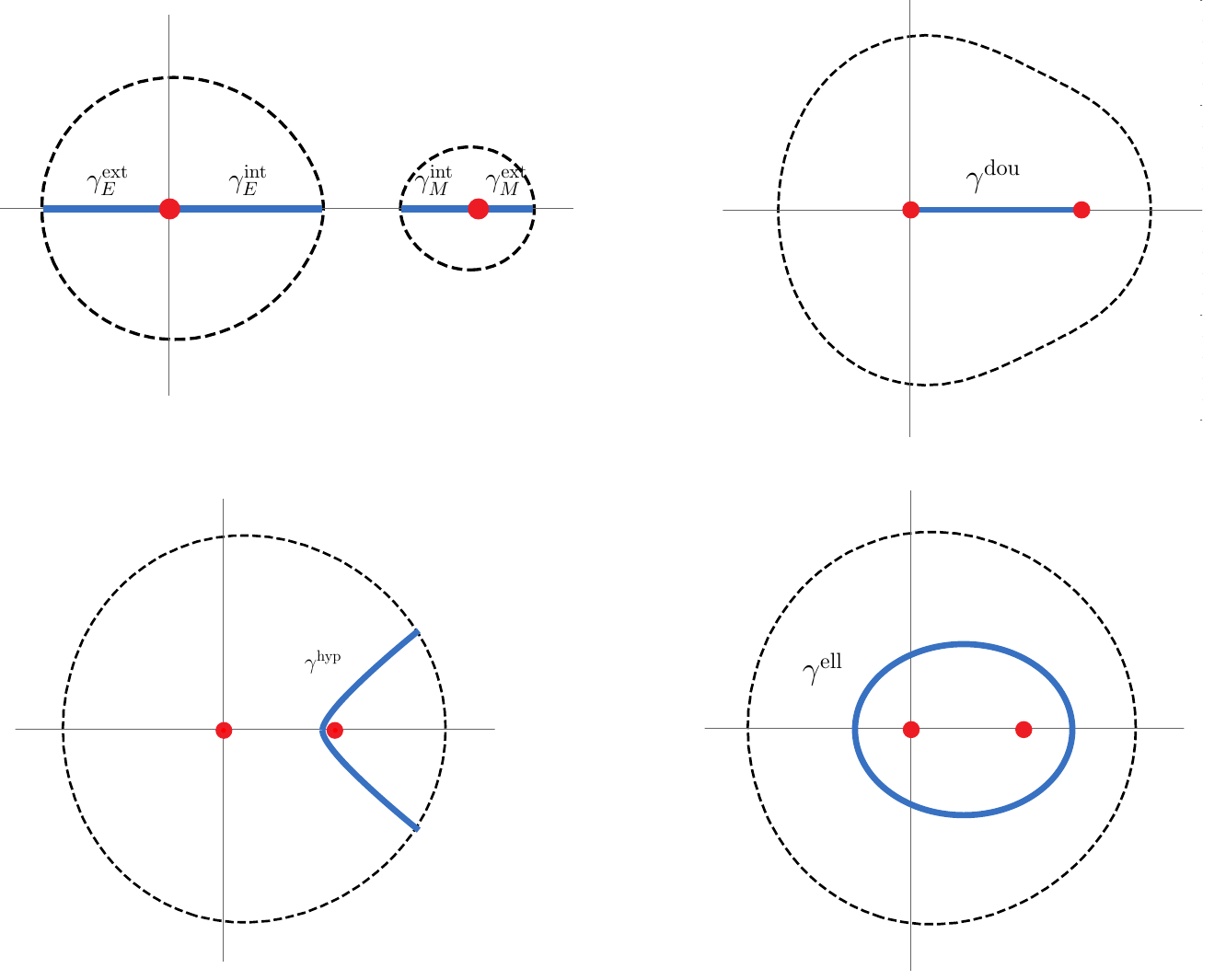}
\caption{Critical orbits in the Euler problem}
\label{cricu}
\end{figure}

For negative energies, there exist four types if $\mu \neq 1/2$ or three types if $\mu =1/2$ of \textit{regular regions}. They are distinguished by five  types  if $\mu \neq 1/2$ or four types if $\mu =1/2$ of \textit{critical curves}, see  Figure \ref{diagram}. Each point on the critical curves represents one of   five types of  \textit{critical orbits}:  the interior and exterior collision orbits near each primary,   the double-collision orbit, the hyperbolic orbit and the elliptic orbit, see Figure \ref{cricu}. The orbits of interest in this paper are the double-collision orbit, the hyperbolic orbit   and the exterior collision orbit   near the Moon: The double-collision orbit exists for $c>c_J$ and up to $c=c_e$ it is stable. However, passing the distinguished energy level $c=c_e$ it becomes unstable. The hyperbolic orbit, which exists for $c \in (c_J, c_h)$, is always unstable. Finally, the exterior collision orbit near the Moon exists for all negative energies and it is unstable if $\mu \neq 1/2$ and  $c\in (c_h, 0)$. In other words, there exist three types of unstable periodic orbits. A natural question then arises: \textit{Do these unstable periodic orbits admit homoclinics?}  Since the Euler problem is integrable, the existence of homoclinics is well-known. In Section \ref{sec:homodd}, we recheck this well-known fact using  the topological invariants of integrable Hamiltonian systems, which will be recalled in  Section \ref{sec:molecule}.

\begin{figure}[h!]
  \centering
  \includegraphics[width=0.76\linewidth]{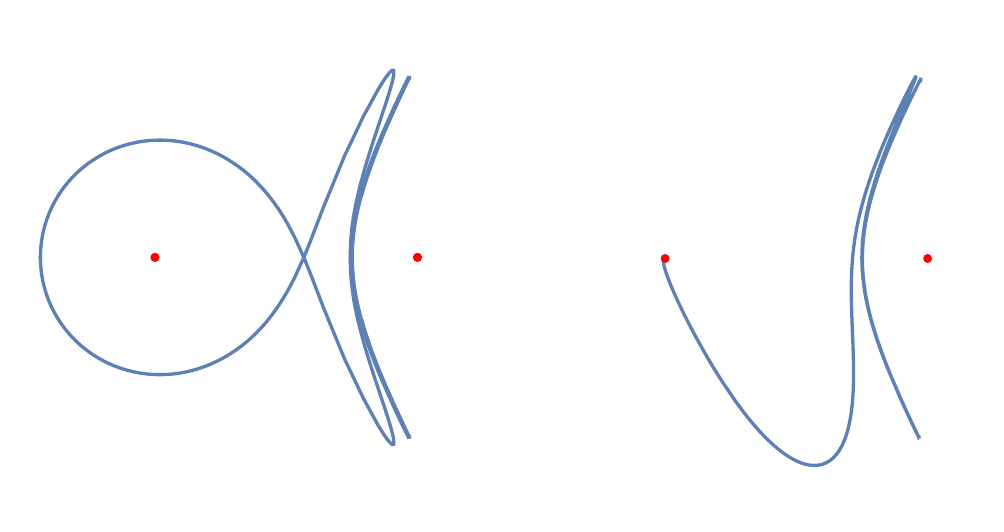}
\caption{Noncollision and collision  homoclinic orbits to the hyperbolic orbit}
\label{homoclnic1}
  \centering
  \includegraphics[width=0.39\linewidth]{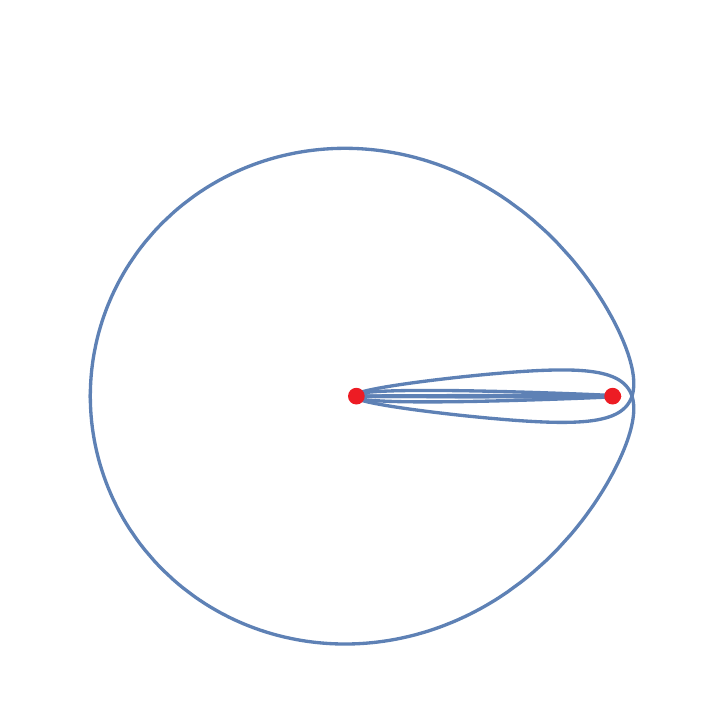}
\caption{A homoclinic orbit to the double-collision orbit}
\label{homoclnicdou}
  \centering
  \includegraphics[width=0.88\linewidth]{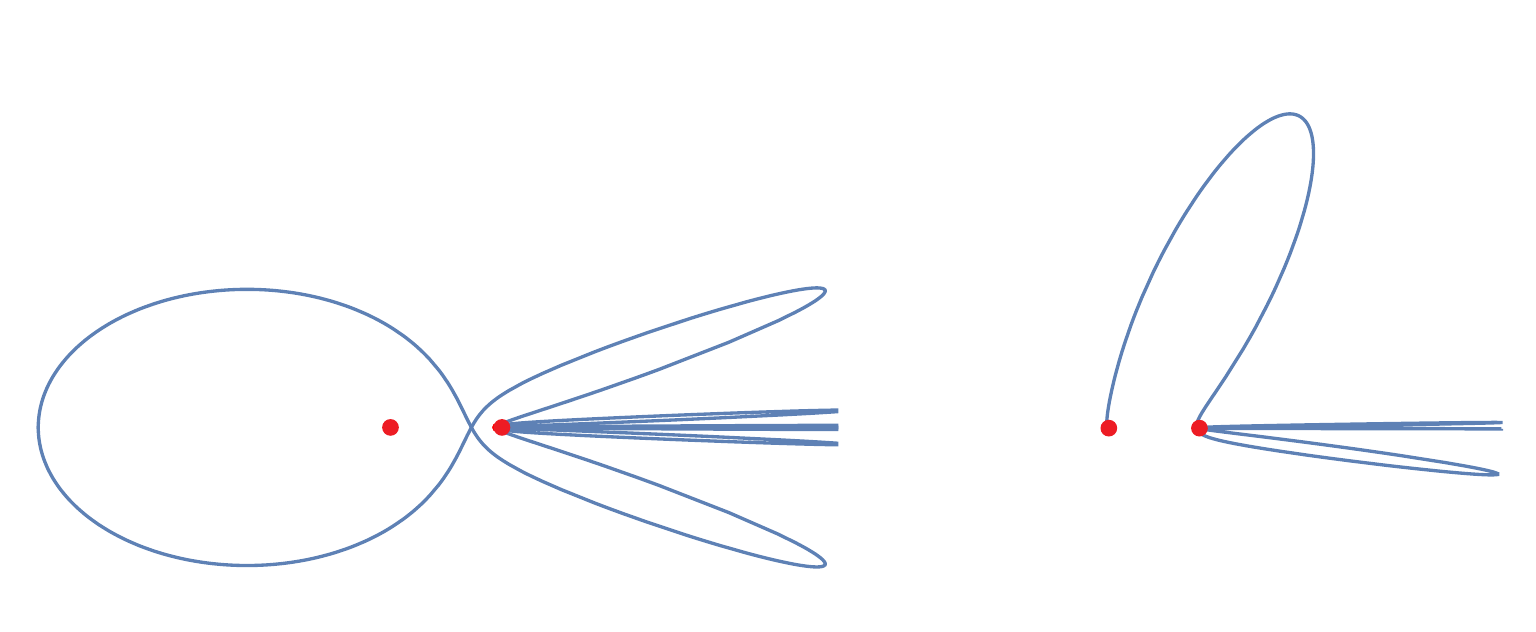}
\caption{Noncollision and collision  homoclinic orbits to the exterior collision orbit in the Moon component}
\label{homoclnicext}
\end{figure}

The main theorem of this paper concerns the behavior of homoclinics to each unstable periodic orbits:

\medskip

\noindent \textbf{Theorem.} The following   do not depend on the choice of the mass ratio:
  \begin{enumerate}[label=(\roman*)]
\item    For each $c \in (c_J, c_h)$, the hyperbolic orbit $\gamma^{\text{hyp}}$ admits homoclinics which collide with either the Earth or Moon. Moreover, all non-collision homoclinic orbits to $\gamma^{\text{hyp}}$  rotate around one of the primaries precisely once, see Figure \ref{homoclnic1}; 
\item    Fix    $c \in (c_e, 0)$ and consider the double-collision orbit $\gamma^{\text{dou}}$.  All homoclinic orbits    rotate around  $\gamma^{\text{dou}}$.  No homoclinics admit collisions, see Figure \ref{homoclnicdou};
\item   Abbreviate by $\gamma^{\text{ext}}_M$ the exterior collision orbit near the Moon for a fixed $c\in (c_h, 0)$ and $\mu \neq 1/2$.  It admits homoclinics which collide with the Earth. No homoclinics can collide with the Moon. Furthermore, any non-collision homoclinic rotates around the Earth precisely once, see Figure \ref{homoclnicext}. 
\end{enumerate}

\bigskip

\noindent \textbf{Remark.} The satellite can rotate   in both direction along all non-collision homoclinic orbits. This  follows from  the fact that the system is invariant under   the anti-symplectic involution $(q,p) \mapsto (q,-p)$. 

 \bigskip

\noindent \textbf{Acknowledgements:}  I wish to express my deepest gratitude to  Urs Frauenfelder for his encouragement. I also thank  Holger Dullin for helpful comments,  the unknown referee for valuable comments and   the Institute for Mathematics of the University of Augsburg for providing a supportive research environment. This work was supported by  Deutsche Forschungsgemeinschaft grants CI 45/8-1 and FR 2637/2-1.

\section{The molecule theory}
\label{sec:molecule}

Let $(M, \omega)$ be a four-dimensional symplectic manifold and $H \in C^{\infty}(M,\R)$ be a Hamiltonian. We assume that the Hamiltonian system $(M, \omega, H)$ is integrable and abbreviate by $F : M \rightarrow \R$ a first integral.   By the Arnold-Liouville theorem, see for example \cite{Arnold}, each compact connected component of  regular common level sets of $H$ and $F$  is a two-dimensional torus, called a \textit{Liouville torus}.  Consequently, the manifold $M$ is foliated by regular leaves which are Liouville tori  and singular leaves along which $dH$ and $dF$ are linearly dependent. This foliation is referred to as the \textit{Liouville foliation}.

\begin{figure}[t]
 \centering
 \includegraphics[width=0.3\textwidth, clip]{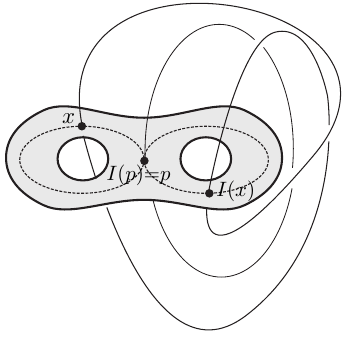}
 \caption{The neighrborhood ${\mathcal{N}}$ of the figure eight equipped with the  involution $I$ }
\label{sudf}
\end{figure}

\begin{figure}[h]
 \centering
 \includegraphics[width=0.8\textwidth, clip]{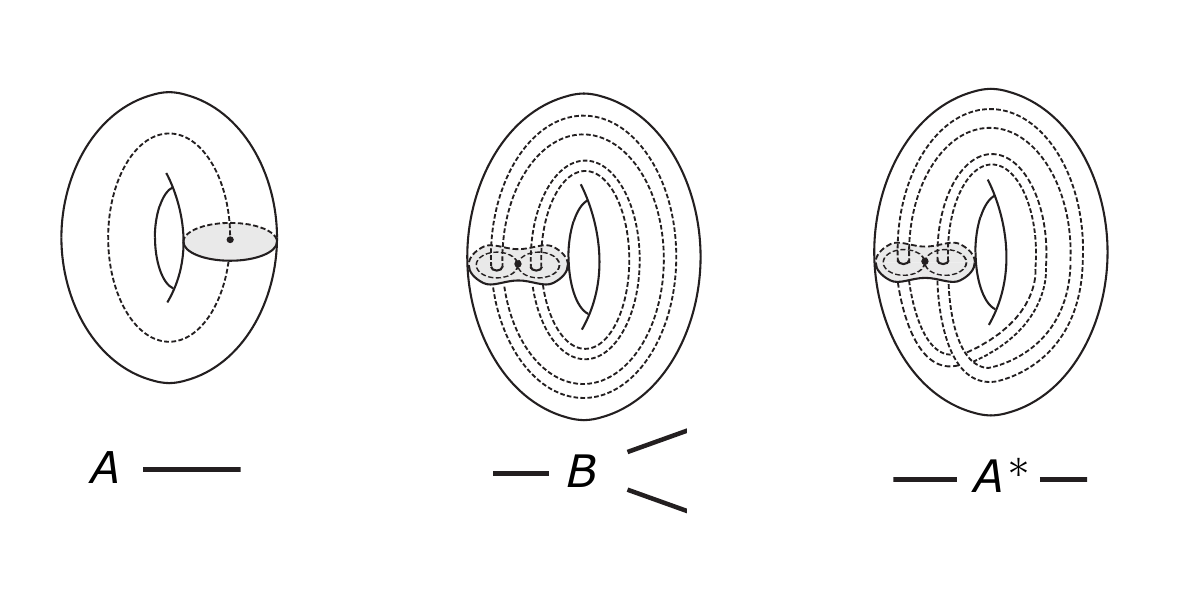}
 \caption{(left) the atom $A$, (middle) the atom $B$, (right) the atom $A^*$}
\label{sudfddd}
\end{figure}

Let $\Sigma$ be a compact regular level set of $H$ which is a three-dimensional submanifold of $M$. By abuse of notation, we use the same symbol $F$ for the restriction $F|_{\Sigma}$. We assume that  $F$ is Morse-Bott, i.e.,  each connected component of the set of critical points of $F$, denoted by $\text{crit} F$,  is a submanifold of $\Sigma$ and   the restriction of $F$ to a small transversal section to $\text{crit} F$ is Morse.    In order to describe the  Liouville   foliation on  $\Sigma$, we recall the Fomenko-Zieschang invariant, see for example \cite{molecule, original}:     Let $L$ be a singular leaf of the Liouville foliation and denote by $U(L)$ its small neighborhood. The neighborhood $U(L)$ is said to be \textit{Liouville equivalent} to another  neighborhood $U'(L)$ of $L$  if there exists a homeomorphism $\phi  : U(L) \rightarrow U'(L)$ which maps the leaves to the leaves. A Liouville equivalence class of  small neighborhoods of $L$ is called an \textit{atom} (associated to $L$).   The number of  singular leaves in an atom is called its {complexity}. In this section, we only consider \textit{simple atoms}, i.e., atoms with complexity one. In view of the classification of atoms, a simple atom has the type either $A$, $B$, or $A^*$, which are described in the following:

\begin{enumerate}[label=(\roman*)]
 \item \textit{the atom  A} is the solid torus $S^1 \times D^2$ whose core $S^1 \times \left \{0 \right\}$ is the singular fiber.  This atom corresponds to either the maximum or   minimum of the Morse-Bott integral $F$ and describes birth or death of a Liouville torus; 
\item \textit{the atom B} is the direct product of a neighborhood $\mathcal{N}$ of the figure eight   and the circle, where the product of the figure eight  and the circle represents the singular fiber. This atom is associated to saddle points of $F$ and describes the decomposition of a Liouville torus into two tori or the reverse; 
\item  \textit{the atom $A^*$} is obtained by glueing the endpoints of the product $\mathcal{N} \times [0,1]$ via the involution $I$ given in Figure  \ref{sudf} in such a way that $(x,0)$ and $( I(x), 1)$ are identified.   This atom is also  associated to saddle points of $F$ and describes  the transition from a Liouville torus into another one, see Figure \ref{sudfddd};
\end{enumerate}

In order to see topology of the three-dimensional manifold $\Sigma$ one needs to know how neighborhoods of   singular fibers are glued along their boundaries.  A combination of atoms of $\Sigma$ is called the \textit{molecule} of $\Sigma$.  Note that a molecule might describe several topologies. For example, if $\Sigma$ has the molecule $A  -A$, then it can be obtained by glueing two solid tori along boundaries and hence it is homeomorphic to either $S^3$, $\R P^3$, $S^1 \times S^2$, or lens spaces. To know the precise topology of $\Sigma$, one needs additional data,  see \cite{molecule, original} and references cited therein.

\section{The Euler problem of two fixed centers}
\label{euler}

In this section we recall some facts on the Euler problem. We introduce   the doubly-covered elliptic coordinates $(\lambda, \nu, p_{\lambda}, p_{\nu})$ which are defined by
\begin{equation*}
\cosh \lambda = |q-E| +|q-M| \in [1,\infty) ,\;\;\;\;\; \cos \nu = |q-E| -|q-M| \in [-1,1],
\end{equation*}
where $(\lambda, \nu) \in \R \times S^1[-\pi, \pi]$. Note that coordinate lines $\lambda = $ constant resp. $\nu=$ constant represent   ellipses resp.   hyperbola in the $q$-plane.  The momenta $p_{\lambda}$ and $p_{\nu}$ are determined by the canonical relation $p_1 dq_1 +p_2 dq_2 = p_{\lambda} d\lambda +p_{\nu} d\nu $ and then the   Hamiltonian becomes 
\begin{equation*}
H = \frac{ H_{\lambda} + H_{\nu}}{\cosh^2 \lambda - \cos^2 \nu},
\end{equation*}
where $H_{\lambda} = 2p_{\lambda}^2 - 2\cosh \lambda$ and $H_{\nu} = 2p_{\nu}^2 + 2(1-2\mu)\cos\nu$.  A Morse-Bott integral is given by
\begin{equation*}
G = - \frac{ H_{\lambda} \cos^2 \nu + H_{\nu} \cosh^2 \lambda}{\cosh^2 \lambda - \cos^2 \nu }. 
\end{equation*}
Given $(G,H)=(g,h)$ the momenta can be expressed by
\begin{equation} \label{newmomentum}
p_{\lambda}^2 = \frac{ c  \cosh^2 \lambda + 2\cosh \lambda + g }{2}, \quad \quad p_{\nu}^2 = \frac{ -c \cos ^2 \nu - 2(1-2\mu)\cos  \nu - g }{2}.
\end{equation}
In particular, the system becomes separable in the elliptic coordinates.

Given $H=c$, we define the new Hamiltonian 
$$
Q = (H-c)(\cosh^2 \lambda- \cos^2 \nu) =  2p_{\lambda}^2 - 2 \cosh\lambda - c \cosh^2 \lambda + 2p_{\nu}^2 + 2(1-2\mu)\cos \nu + c \cos^2 \nu.
$$
Note that orbits of $H$ with energy $c$ and the time parameter $t$ correspond to orbits of $Q$ with energy $0$ and the time parameter $\tau$, which is defined by
$$ \tau = \int \frac{ dt}{\cosh^2 \lambda -\cos^2\nu}.$$
Moreover,   each bounded component of the energy hypersurface $H^{-1}(c)$ compactifies to each compact component, which is diffeomorphic to $S^3$, of the regularized hypersurface $Q^{-1}(0)$.   In particular, the satellite is allowed to pass through the primaries. In the following we assume that the system is regularized.

We now define the function  
$$f_{\mu}(x) = cx^2 + 2(1-2\mu)x+g$$
 so that $p_{\lambda} = \pm \sqrt{  { f_0(\cosh \lambda)}/{2}}$ and $ p_{\nu} = \pm \sqrt{-  { f_{\mu}(\cos \nu)}/{2}}. $   
The classically allowed values $(g,c)$ are those at which the momenta are real. Note that for the momenta to be real we need $f_0 >0$ for $ x \in [1, \infty)$ and $f_{\mu}<0$ for $\mu \neq 0 $ and $x \in [-1, 1]$.  Abbreviate by $x_1^{\mu}$, $x_2^{\mu}$ two roots of $f_{\mu}$. 
   In the following we describe the procedure   to find all such values in the lower half $(g,c)$-plane, which is given for example in \cite{Strand, Bifurcation}.

\;\;

\textit{Case 1.} the function $f_{\mu}$ has no real roots.\\
In this case we have $(1-2\mu)^2 < gc$ and $f_{\mu} <0$. 

\;\;

\textit{Case 2.} the function $f_{\mu}$ admits a common root $x^{\mu}$.\\
This implies that $(1-2\mu)^2 =gc$ and $x= - ({1-2\mu})/{c}$.

 \;\;

\textit{Case 3.} the function $f_{\mu}$ has two real roots $ x_1^{\mu} < x_2^{\mu}$. In the following we omit the obvious condition $(1-2\mu)^2 >gc$.
\begin{enumerate}[label=(\roman*)]
\item $ x_1^{\mu} < x_2^{\mu} < -1 < 1$: since $x_2^{\mu} = {( -(1-2\mu)- \sqrt{(1-2\mu)^2 -gc})}/{c}$ is positive, this is impossible;
\item $x_1^{\mu} < -1 < x_2^{\mu}< 1 \; \Leftrightarrow \; 2(1-2\mu) < c+g<-2(1-2\mu)$ which is impossible;
\item $x_1^{\mu} < -1 < 1< x_2^{\mu} \; \Leftrightarrow \; 0< c+g-2(1-2\mu)$;
\item $-1<  x_1^{\mu}  < x_2^{\mu}< 1\; \Leftrightarrow \;  c+g+2(1-2\mu)<0$ with $c<-(1-2\mu)$; 
\item $-1 < x_1^{\mu} < 1 < x_2^{\mu}\; \Leftrightarrow \; -2(1-2\mu)< c+g< 2(1-2\mu)$;
\item $-1< 1 < x_1^{\mu}  < x_2^{\mu}\; \Leftrightarrow \; c+g+2(1-2\mu)<0$ with $c>-(1-2\mu)$;
\item $x_1^{\mu}=-1 < x_2^{\mu} < 1 \; \Leftrightarrow \;   4(1-2\mu)<0$   which is impossible;
\item $x_1^{\mu}=-1 < x_2^{\mu} = 1 \; \Leftrightarrow \; c+g=\pm2(1-2\mu)$ which is impossible;
\item $x_1^{\mu}=-1 < 1 <x_2^{\mu}  \; \Leftrightarrow \; c+g-2(1-2\mu)=0$ if $\mu \neq {1}/{2}$;
\item $x_1^{\mu}<-1 < x_2^{\mu} = 1 \; \Leftrightarrow \; 0<- 4(1-2\mu)$  which is impossible;
\item $-1 < x_1^{\mu} < x_2^{\mu} = 1 \; \Leftrightarrow \;  c+g+2(1-2\mu)=0$ with $c<-(1-2\mu)$;
\item $-1 < 1 = x_1^{\mu} < x_2^{\mu} \; \Leftrightarrow \; c+g+2(1-2\mu)=0$ with $c>-(1-2\mu)$.
\end{enumerate}

\;\;

In view of the fact that for the momenta to be real we need $f_0>0$ with $x \in [1,\infty)$ and $f_{\mu}<0$, $\mu \neq 0$, with $x \in [-1,1]$, we combine the above results and give the bifurcation diagram in Table \ref{tablerange} and Figure \ref{diagram}.
\begin{table}[h]
\begin{tabular*}{\textwidth}{c    p{6.5cm}    c}
\toprule

  \textbf{\small{Region}} &  \centering \textbf{\small{Ranges of the Roots}}&   \textbf{\small{Ranges of the Variables}}   \\

\midrule \\

$P$ &  \footnotesize{\makecell{ \noindent $\cosh\lambda : -1 < 1 < x_1^0 <x_2^0$ \hspace{2.8cm} \\ $ \cos\nu : \begin{cases} -1 < 1 < x_1^{\mu} < x_2^{\mu} ,   &  (1-2\mu)^2 \geq gc  \\  \hfil \text{complex roots} ,&   (1-2\mu)^2 < gc   \end{cases}$ }}& \footnotesize{$ \begin{matrix}  \cosh \lambda \in   [ x_1^0, x_2^0 ] \\  \\ \cos \nu \in    [-1, 1]  \end{matrix}$}  \\

\\

\\

$L$ &   \footnotesize{\makecell{ \noindent $\cosh\lambda : -1 <  x_1^0 < 1 < x_2^0$ \hspace{2.8cm} \\ $ \cos\nu : \begin{cases} -1 < 1 < x_1^{\mu} < x_2^{\mu} ,   &   (1-2\mu)^2 \geq gc  \\  \hfil \text{complex roots}, &   (1-2\mu)^2 < gc   \end{cases}$              }}&      \footnotesize{  $   \begin{matrix}  \cosh \lambda  \in   [ 1, x_2^0 ] \\  \\   \cos \nu \in [-1, 1]  \end{matrix} $}  \\

\\

\\

$S$ &  \footnotesize{\makecell{ \noindent $\cosh\lambda : -1 <  x_1^0< 1 <x_2^0$\hspace{2.8cm} \\ $ \cos\nu :   -1 < x_1^{\mu}   < x_2^{\mu} <1   $     \hspace{2.8cm}         }}&     \footnotesize{   $  \begin{matrix}  \cosh \lambda   \in   [ 1, x_2^0 ] \\  \\  \cos \nu \in  ( [-1, x_1^{\mu}] \cup [x_2^{\mu}, 1] )   \end{matrix} $}  \\

\\

\\

$S'$ & \footnotesize{ \makecell{ \noindent $\cosh\lambda : -1  < x_1^0 < 1<x_2^0$ \hspace{2.7cm} \\ $ \cos\nu :    -1 < x_1^{\mu} < 1   < x_2^{\mu}        $ \hspace{2.8cm}    }}&   \footnotesize{    $  \begin{matrix}  \cosh \lambda \in   [ 1, x_2^0 ] \\  \\   \cos \nu \in [-1, x_1^{\mu}]   \end{matrix} $ } \\

 \\

\bottomrule
\end{tabular*}

\caption{The ranges of the roots and the variables in the four regular regions}

\label{tablerange}
\end{table}
 
As mentioned in the introduction, if $\mu \neq 1/2$,  in  the lower-half $(G,H)=(g,c)$-plane there exist   four  {regular regions} consisting of regular values of the   energy momentum mapping $(\lambda, \nu) \mapsto (G(\lambda, \nu), H(\lambda, \nu))$. By the Arnold-Liouville theorem, each point in the regular regions represents a Liouville torus. Following the notations from \cite{Charlier, Pauli}, the regular regions are labeled by ${S'}$, ${S}$(satellite), ${L}$(lemniscate), and ${P}$(planetary).  They are bounded by the following five \textit{critical curves}:
\begin{eqnarray*}
&&\ell_{1,2} : c = -g \pm 2(1-2\mu), \;  \quad \quad \quad \quad \quad \quad  \ell_3 : c=-g-2,\\
&& \ell_4 : gc=(1-2\mu)^2,~c_J<c<c_h,  \quad \quad \quad \ell_5 : gc=1,~c_e <c.
\end{eqnarray*}
Here,  $c_e=-1$ or  $c_h=-(1-2\mu)$  are the energy levels at which the line $ \ell_3$ and   curve $ \ell_5$ or the line $ \ell_2$ and   curve $ \ell_4$ intersect, respectively(the letters $e$ and $h$ stand for \textit{elliptic} and \textit{hyperbolic}). Note that at these points the Liouville foliation changes.     If the two primaries have the equal masses, $\mu=1/2$, then the two curves $ \ell_1$, $ \ell_2$ become equal  and hence the $S'$-region does not appear.  In the $S$-region each point represents a motion of the satellite which is confined to a neighborhood of either the  Earth or Moon,  while in the $S'$-region the  motion takes place only near the Earth, see for example   \cite{Strand, Bifurcation}. 

We now examine the critical orbits:
\begin{enumerate}[label=(\roman*)]
\item  \textit{On $l_1 : c = -g+2(1-2\mu)$}. The satellite moves along the ray $\cos\nu = -1$, but the motion is bounded by the ellipse $\cosh\lambda = x_2^0$. We call this orbit the \textit{exterior collision orbit} in the Earth component;  
\item \textit{On $l_2 : c = -g-2(1-2\mu)$}. Each point on the line $l_2$ represents an orbit in a neighborhood of either the Earth or  Moon. Assume the satellite moves near the Earth so that $(\cosh\lambda, \cos\nu) \in[1, x_2^0] \times [-1, x_1^{\mu}]$. In particular, the motion is regular. Consider orbits near the Moon. Then the satellite moves along the ray $\cos\nu = 1$, but bounded by the ellipse $\cosh\lambda = x_2^0$. This orbit is called the \textit{exterior collision orbit} in the Moon component;
\begin{figure}[t]
  \centering
  \includegraphics[width=0.6\linewidth]{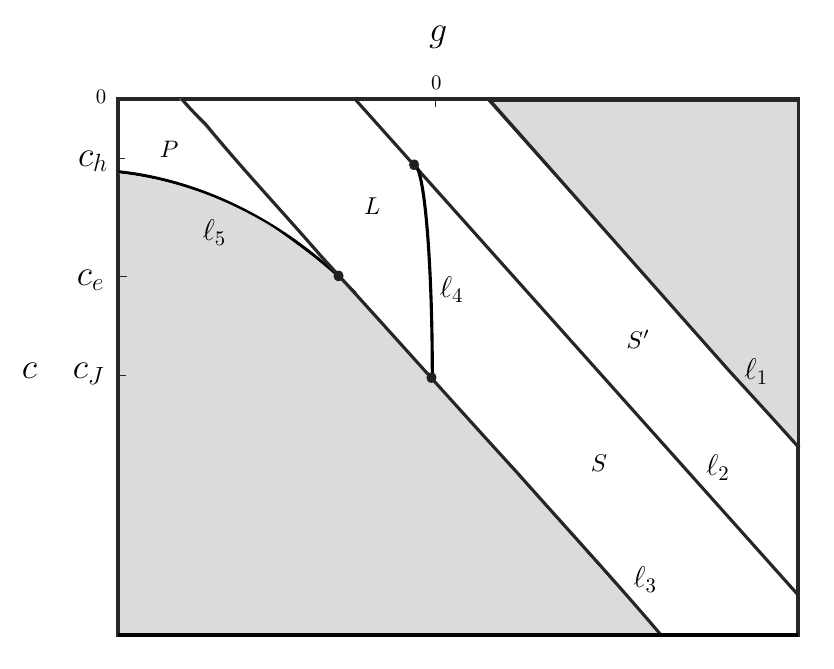}
\caption{The case $\mu\neq 1/2$: The lower half $(g,c)$-plane is divided into four regions by the five critical curves. The shaded regions are classically not allowed.  For negative energies, each regular level of the energy-momentum mapping represents one of  four types of motions, which are labeled by $S'$, $S$, $L$ and $P$.}
\label{diagram}
\end{figure}
\item \textit{On $l_3: c = -g-2$}. On the line $l_3$ we have $\cosh\lambda=1$, i.e, the satellite moves along the line segment joining the two masses. If $c<c_J$, the variable $\cos\nu$ lies in $[-1, x_1^{\mu}] \cup [x_2^{\mu}, 1]$, where the first or the second interval corresponds to the motion  in the Earth or  Moon component, respectively. For $\cos\nu \in[-1, x_1^{\mu}]$, the satellite moves on the line segment $\left \{ (q_1, 0) :  -1/{2} \leq q_1 \leq b_c^E \right \}$, where $(b_c^E, 0) \in \partial \mathcal{K}_c^E$ with $b_c^E \in (-1/{2},l)$.  Similarly, in the Moon component we have $\left \{ (q_1, 0) : b_c^M \leq q_1 \leq 1/2 \right \}$, where $(b_c^M, 0) \in \partial \mathcal{K}_c^M$ with $b_c^M \in (l,{1}/{2})$. We refer to such an orbit  as the \textit{interior collision orbit} in the Earth or  Moon component. As $c$ passes over $c_J$, the two components of the Hill's region become connected. The two interior collision orbits also become connected so that the satellite can move between the Earth and Moon. Indeed, there is no restriction on $\cos\nu$, i.e., $\cos\nu \in [-1,1]$ for $c>c_J$.  This orbit is referred to as the \textit{double-collision orbit};
\item  \textit{On $l_4 : gc =  (1-2\mu)^2$}. We have $\cos\nu \in [-1, - ({1-2\mu})/{c})$, $\cos\nu \in (- ({1-2\mu})/{c}, 1]$, or $\cos\nu = - ({1-2\mu})/{c}$.   For the last  case, the satellite moves along the hyperbola $\cos\nu = - ({1-2\mu})/{c}$ within the ellipse $\cosh\lambda = x_-^0$. We call this orbit the \textit{hyperbolic orbit}. 
Each member of the family $\gamma_{\text{hyp}}^c$  of the hyperbolic orbits, $c \in (c_J, c_h)$,  is the Lyapunov orbit. To see this, we need to show that the family $\gamma_{\text{hyp}}^c$ converges uniformly to $ L = (l,0,0,0 )$ as $c$ tends to $c_J$ from above, i.e, it emanates from the critical point. Recall that along the hyperbolic orbits  the equation $c\cos\nu^2 + 2(1-2\mu) \cos\nu +g =0$ has double roots $\cos \nu = - ({1-2\mu})/{c}$ which implies that   the hyperbola  is given by
\begin{equation*}
 |q+(\frac{1}{2},0)| - |q - (\frac{1}{2},0)| = - \frac{1-2\mu}{c}.
\end{equation*}
This hyperbola, for $\mu \neq 1/2$, is closer to the Moon than to the Earth since $-({1-2\mu})/{c}>0$. It remains to show that     the hyperbola $\cos\nu= -({1-2\mu})/{c}$ consists of the single  point $(l,0)$ at $c=c_J$, or equivalently
\begin{eqnarray*}
\label{set}
\left \{ (q_1, q_2)  : |q+(\frac{1}{2},0)| - |q-(\frac{1}{2},0)| = - \frac{1-2\mu}{c_J} \right \} \cap \mathcal{K}_{c_J} = \left \{ (l,0) \right \}.
\end{eqnarray*}
Since $\cosh\lambda \rightarrow 1 $ as $c \rightarrow c_J$, it suffices to show that the vertex of this hyperbola  is given by $(l,0)$. We compute that

\begin{equation*}
  \frac{1-2\mu}{ -c_J} = \frac{ 1-2\mu}{ 1+2\sqrt{\mu(1-\mu)} } = \frac{ 1 -2\sqrt{\mu(1-\mu)} }{1-2\mu}.
\end{equation*}
Then the vertex of this hyperbola   is given by
\begin{equation*}
(q_1, q_2) = \bigg( \frac{1-2\sqrt{\mu(1-\mu)}}{2(1-2\mu)}, 0 \bigg) = (l,0).
\end{equation*}
This completes the proof of the claim.  On the other hand, as $c \rightarrow c_h$ from below the hyperbolic orbits degenerates to the exterior collision orbit in the Moon component;
\item  \textit{On $l_5 : gc = 1$}.  We have $\cosh\lambda = - {1}/{c}$, but no restriction on $\cos\nu$. This implies that the satellite moves along the ellipse $\cosh\lambda = - {1}/{c}$. This orbit is referred to as the \textit{elliptic orbit}. As $c \rightarrow c_e$ from above,  the ellipse $\cosh\lambda=-{1}/{c}$  degenerates to the line segment $\cosh\lambda=1$ joining the two masses, namely the double-collision orbit.
\end{enumerate}

\section{Molecules of the Euler problem}
\label{sec:eulermolecule}

\begin{figure}[h]
 \centering
 \includegraphics[width=0.8\textwidth, clip]{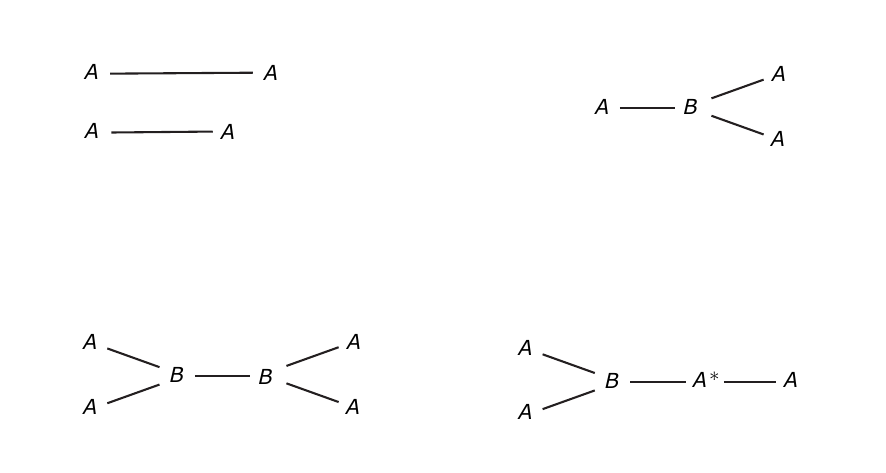}
 \caption{ Molecules of the Euler problem. (left above) For $c<c_J$:  The above one is associated to the Earth component and the below one is associated to the Moon component. The leftmost $A$-atoms  represent the interior collision orbits  and the rightmost ones represent the exterior collision orbits. (right above) For $c_J<c_e$: The leftmost $A$ represents the double-collision orbit and the middle $B$ is the hyperbolic orbit. The right two atoms are exterior collision orbits in both components.  (left below) For $c_e < c < c_h$: The leftmost atoms are the elliptic orbits. The first, from the left, atom $B$ is the double-collision orbit and the second one is the hyperbolic orbit. The rightmost atoms are the exterior collision orbits. (right below) For $c_h<c<0$: As in the previous case, the first two atoms are the elliptic orbits and the atom $B$ is the double collision orbit. In this energy range, the hyperbolic orbit does not appear and the exterior collision orbit in the Moon component becomes unstable: the atom $A^*$. The last atom $A$ represents the exterior collision orbit in the Earth component.   }
\label{tdde}
\end{figure}

The discussion in the previous section together with the results of     the molecule theory described in Section \ref{sec:molecule} give rise to the molecules of the Euler problem, which are already  given by Waalkens-Dullin-Richter in \cite[Section 2]{Bifurcation}, see Figure \ref{tdde}.  We refer to \cite{sphere} for the molecule structure of the Euler problem on the two-sphere.

\section{The existence of homoclinic orbits}
\label{sec:homodd} 
 Recall that in the Euler problem the hyperbolic orbit $\gamma^{\text{hyp}}$ for $c\in (c_J, c_h)$, the double-collision orbit $\gamma^{\text{dou}}$ for $c \in (c_e ,0)$ and the exterior collision orbit in the Moon component $\gamma^{\text{ext}}_M$ for $c \in (c_h, 0)$, if $\mu \neq 1/2$,  are unstable periodic orbits. The existence of homoclinic orbits to each of these   orbits follows immediately from the properties of the corresponding atoms as follows:

\medskip

\textit{Case 1. The hyperbolic orbit $\gamma^{\text{hyp}}$}\\
Recall that the corresponding atom is the type $B$ whose singular fiber, denoted by $L$,  is homeomorphic  to the union of two two-dimensional tori attached along a single periodic orbit which is in this case the hyperbolic orbit $\gamma^{\text{hyp}}$. A manifold $L \setminus \gamma^{\text{hyp}}$  consists of two connected components, where each of them is diffeomorphic to the two-torus minus a periodic orbit:   one is around the Earth and the other is around the Moon.  Here we  consider only  the component near the Earth, which is abbreviated by $L_E \cong  \Bbb{T}^2 \setminus \gamma^{\text{hyp}}$. For the component near the Moon, the same argument also holds. Since $\gamma^{\text{hyp}}$ is periodic, it has some rational slope on the torus. We then observe that on $L_E$ there exist no periodic orbits. Otherwise, all periodic orbits have the same rational slope with $\gamma^{\text{hyp}}$, but this is not the case in view of the symmetry associated with the integral $G$. Therefore, every orbit on $L_E$ is of irrational slope and hence is dense. This implies that every orbit converges asymptotically to $\gamma^{\text{hyp}}$ in  both directions, namely all   orbits on $L_E$ are homoclinic orbits.

\medskip

\textit{Case 2. The double-collision orbit $\gamma^{\text{dou}}$}\\
In this case we also have the atom $B$. The argument given in the previous case carries over to this case and hence there exist homoclinic orbits to $\gamma^{\text{dou}}$.

\medskip

\textit{Case 3. The exterior collision orbit in the Moon component $\gamma^{\text{ext}}_M$}\\
We have the atom $A^*$. The singular fiber $L$ is given by the "skew product" of the figure eight and the circle via the involution given in Figure \ref{sudf}.  Note that  $L \setminus \gamma^{\text{ext}}_M$ consists of a single connected component: $L \setminus \gamma^{\text{ext}}_M \cong \Bbb{T}^2 \setminus \gamma^{\text{hyp}}$. The remaining argument is the same as before.

\begin{Remark} \rm The above discussions on the existence of homoclinic orbits   only  make use of    the existence of the atoms $B$ and $A^*$. Therefore, for integrable Hamiltonian systems with two degrees of freedom, any unstable periodic orbit  associated to the atom $B$ or $A^*$  admits homoclinic orbits.
\end{Remark}

\section{Proof of the theorem}
\label{sec:homo}

Since homoclinic orbits to a unstable periodic orbit $\gamma$  lie on the intersection of the unstable and   stable manifolds of $\gamma$, if they exist, then   they take  the same integral value  as $\gamma$. 

\medskip

\textit{Case 1. The hyperbolic orbit $\gamma^{\text{hyp}}$}\\
 Let us fix the energy level $c \in (c_J, c_h)$ at which $\gamma^{\text{hyp}}$ exists. Recall from  Section \ref{euler} that the corresponding value $(g,c)$ lies on $\ell_4$ so that  $\gamma^{\text{hyp}}$ takes $g_c ^{\text{hyp}}:=  {(1-2\mu)^2}/{c} <0$. 
 
We first consider homoclinic orbits along which collisions do not occur and observe the behavior of the satellite along these orbits.  In view of the argument in Section \ref{euler},  such orbits lie in a neighborhood of the Earth or  Moon in the case $\cos\nu \in [-1, -({1-2\mu})/{c})$ or $\cos\nu \in (  -({1-2\mu})/{c} , 1]$, respectively, where   $\cos\nu = -({1-2\mu})/{c}$ represents $\gamma^{\text{hyp}}$. In $(\lambda, \nu)$-coordinates, the range for $\nu$ is given by $[-\pi, -\nu_1 ) \cup ( \nu_1 , \pi]$, where $ \cos \nu_1 = -({1-2\mu})/{c}$.  Without loss of generality we only consider the interval $[-\pi, -\nu_1)$. For $( \nu_1 , \pi]$, the satellite shows the same behavior.

By the Hamiltonian equations and \eqref{newmomentum} we have
\begin{eqnarray}
\nonumber\dot{\nu} &=&  \frac{\partial H}{\partial p_{\nu}} \\
\nonumber &=& \frac{4 p_{\nu}}{\cosh^2 \lambda  - \cos^2 \nu}\\
\nonumber&=& \pm \frac{4}{\cosh^2 \lambda - \cos^2 \nu} \sqrt{\frac{  c^2 \cos ^2 \nu + 2(1-2\mu)c \cos  \nu +(1-2\mu)^2 }{-2c} } \\
\label{sdf33dss2} &=& \pm \frac{2\sqrt{2}\sqrt{-c}}{\cosh^2 \lambda - \cos^2 \nu} \bigg( \cos \nu + \frac{1-2\mu}{c} \bigg) .
\end{eqnarray}
Since we are considering non-collision orbits, the denominator $\cosh^2 \lambda- \cos^2 \nu$ never vanishes. Therefore, the sign of $\dot{\nu}$ remains unchanged unless $\nu=\pm \nu_1$. Thus, we conclude that $\dot{\nu}$ is nonvanishing. We also have 
$$
\dot{\lambda} = \pm \frac{4}{\cosh^2 \lambda - \cos^2\nu}\sqrt{\frac{ c^2 \cosh^2 \lambda + 2c \cosh\lambda +(1-2\mu)^2}{2c}}.
$$
Note that $c^2 \cosh^2 \lambda + 2c \cosh\lambda +(1-2\mu)^2$ is negative for $\cosh\lambda \in [-1, c_J/c]$ which shows that $\dot{\lambda}$ is also nonvanishing. Therefore, $\cosh\lambda$ oscillates in the interval $[-1, c_J/c]$.

 Suppose that $\nu = \nu_0 \in (-\pi, -\nu_1)$ and $\dot{\nu} >0$ at $t = t_0$. Then by the previous argument, $\nu $ converges asymptotically to $\nu=-\nu_1$ as $t \rightarrow \infty$. On the other hand, in backward time, i.e., as $t$ decreases, $\nu$ also converges asymptotically to $\nu=-\nu_1$. Since $\nu$ is defined on the circle $S^1[-\pi, \pi]$ this picture holds true.
 We obtain a similar picture if we start with $\dot{\nu} <0$ at $t = t_0$.

Along each non-collision homoclinic orbit  the variable $\nu$ attains $-\pi$ precisely once. In $(q_1, q_2)$-coordinates, this implies that the satellite crosses the subset $K$ of the negative $q_1$-axis, which is the line segment joining the boundary of the Hill's region $\mathcal{K}_c^E$ and the Earth, precisely once. To prove the assertion that every non-collision homoclinic orbit rotates around one of the two primaries precisely once, it remains to show that the trajectory cannot be tangent to $K$.  Suppose that   the tangency occurs at $t = t'$ and  hence $\dot{\nu} = 0$.  By the equation (\ref{sdf33dss2}), that $\dot{\nu}=0$ is equivalent to that $c = c_h = -1+2\mu$, which contradicts to $c \in (c_J, c_h)$, see Figure \ref{homoclnic1}.

To show the existence of collision homoclinic orbits, we need to show that its projection to the configuration space contains the two primaries at which we have $\cosh \lambda = 1$, $\cos \nu = \pm 1$. It then suffices to show that for a fixed  $g=  {(1-2\mu)^2}/{c}$, $ c \in (c_J, c_h)$ and $(\cosh \lambda , \cos\nu ) = ( 1 , \pm 1)$, the squared momenta \eqref{newmomentum} are positive. To see this, we first compute that
\begin{eqnarray*}
p_{\lambda}^2 &=& \frac{c \cosh ^2 \lambda + 2 \cosh\lambda + g}{2} \\
&=& \frac{ c + 2 +  (1-2\mu)^2 /c }{2} \\
&=& \frac{c^2 + 2 c + (1-2\mu)^2 }{2c}  .
\end{eqnarray*}
Since $ c_h = -1+2 \mu \leq -1 + 2 \sqrt{\mu(1-\mu)}$ in view of $\mu \leq 1/2$, we obtain that $c^2 + 2 c + (1-2\mu)^2$ is negative for $ c\in (c_J, c_h)$ which implies that $p_{\lambda}^2 $ is positive. We also compute that
\begin{eqnarray*}
p_{\nu}^2 &=& \frac{ -c \cos ^2 \nu - 2(1-2\mu)\cos  \nu - g }{2} \\
&=& \frac{ -c \pm 2 (1-2\mu) - (1-2\mu)^2 /c}{2} \\
&=& \frac{ ( c \pm (1-2\mu) )^2}{-2c} >0.  
\end{eqnarray*}
This proves the assertion on the existence of collision homoclinic orbits, see Figure \ref{homoclnic1}.

\medskip

\textit{Case 2. The double-collision orbit $\gamma^{\text{dou}}$}\\
We fix $c \in (c_e, 0)$. Recall that  $\gamma^{\text{dou}}$ takes the value $g_c^{\text{dou}}:=-c-2 \in (-2,-1)$.   Note that 
$$
(\lambda, \nu) \mapsto (q_1, q_2) = \bigg( \frac{1}{2} \cosh \lambda \cos \nu, \frac{1}{2} \sinh\lambda \sin \nu \bigg)
$$
is a 2-to-1 covering with two branch points $E$ and $M$. The two sheets are related by the involution $(\lambda , \nu) \mapsto (-\lambda, -\nu)$ which extends to the phase space by 
\begin{equation}\label{eqinvo}
(\lambda, \nu , p_{\lambda}, p_{\nu}) \mapsto (-\lambda, -\nu , -p_{\lambda}, -p_{\nu}). 
\end{equation}
Any homoclinic orbit to $\gamma^{\text{dou}}$ has $\cosh \lambda \in (1, \cosh\lambda_1]$, where $\lambda_1 >0$ satisfies $\cosh \lambda_1 = ( -2-c)/c $. Hence we have $\lambda \in [-\lambda_1, 0) \cup (0, \lambda_1]$.  As a result, the two connected components  of the singular fiber minus $\gamma^{\text{dou}}$ are related to each other via the involution \eqref{eqinvo}.  Therefore, to observe the behavior of the satellite along homoclinic orbits, without loss of generality, we may consider only $\lambda \in (0, \lambda_1]$. 

Following  an argument similar with the one given  in the previous case one can easily see  that there exist no collision homoclinic orbits to $\gamma^{\text{dou}}$. In view of the  Hamiltonian equations and \eqref{newmomentum} along homoclinics we have
$$
\dot{\lambda} = \pm \frac{4}{\cosh^2 \lambda - \cos^2 \nu} \sqrt{   \frac{c \cosh^2 \lambda +2 \cosh\lambda -c-2}{2}     }
$$
which shows that $\dot{\lambda} =0$ if and only if $\lambda = \lambda_1$.  Suppose that $\lambda   \in (1, \lambda_1]$ and $\dot{\lambda} >0$ at $t=t_0$. In forward time,  $\dot{\lambda}$ remains to be positive until $\lambda=\lambda_1$ so that $\lambda$ increases. At $\lambda=\lambda_1$ the velocity $\dot{\lambda}$ vanishes and then becomes negative. As time further increases, $\lambda$ decreases and converges asymptotically to $\lambda=1$. In backward time $\lambda$ also converges asymptotically to $\lambda=1$. We obtain a similar picture for the other case.

Note that  there are no constraints on the variable $\nu$, i.e.,  $\nu \in S^1[-\pi, \pi]$. As in the previous case, one can easily show that   $\dot{\nu}$ is nonvanishing along homoclinics.   As a result, $\dot{\nu}$ is either always positive or always negative which shows that  along  homoclinic orbits  $\cos \nu$ oscillates in the interval $[-1,1]$. This together with the previous discussion show that homoclinic orbits rotate around the double-collision orbit, see Figure \ref{homoclnicdou}.

\medskip

\textit{Case 3. The exterior collision orbit in the Moon component $\gamma^{\text{ext}}_M$}\\
We fix $\mu <1/2$ and $c \in (c_h, 0)$. Abbreviate $g_c^{\text{ext}}= -c-2(1-2\mu)$. In this case we have $\cos \nu \in [-1,1)$, where $\cos \nu =1$ represents $\gamma^{\text{ext}}_M$. The range for $\nu$ is then given by $[-\pi, 0) \cup (0, \pi]$. Without loss of generality, we only consider $[-\pi, 0)$. 

By the same reasoning as in the first case with $\nu_1=0$, we obtain that $\dot{\lambda}$ and $\dot{\nu}$ are nonvanishing along noncollision homocilinics and hence $\lambda$ oscillates and $\nu$ converges asymptotically to $0$ in backward and forward time. Moreover, each non-collision homoclinic crosses the line segment,  which is a subset of the negative $q_1$-axis and which joins the Earth and the boundary $ \partial \mathcal{K}_c$,  precisely once. We also obtain that there exist   homoclinics which collide with the Earth. There exist no collisions with the Moon, see  Figure \ref{homoclnicext}.

\;\;\;

\end{document}